\documentclass[12pt]{article}

\usepackage{amsfonts,amssymb,amsmath,amsthm,eucal}
\usepackage[dvips]{graphics}

\numberwithin{equation}{section}

\newcommand{\R}{\mathbb{R}}
\newcommand{\C}{\mathbb{C}}

\newcommand{\GG}{\mathsf{G}}
\newcommand{\Gc}{\mathsf{G}^\circ}

\newcommand{\HH}{\mathsf{H}}
\newcommand{\KK}{\mathsf{K}}

\newcommand{\RR}{\mathsf{R}}
\newcommand{\Rc}{\mathsf{R}^\circ}
\newcommand{\EE}{\mathsf{E}}

\newcommand{\erf}{\operatorname{erf}}

\newcommand{\bF}{\mathbf{F}}

\newcommand{\cE}{\mathcal{E}}
\newcommand{\cEt}{\cE^\circlearrowleft}

\newcommand{\cF}{\mathcal{F}}
\newcommand{\cG}{\mathcal{G}}
\newcommand{\ad}{\dot a}
\newcommand{\Ad}{\dot A}
\newcommand{\e}{\varepsilon}

\newcommand{\af}{\mathsf{a}}
\newcommand{\A}{\mathsf{A}}
\newcommand{\wtt}{\widetilde{\tau}}

\newcommand{\ix}{x}

\newtheorem{theorem}{Theorem}
\newtheorem{lemma}{Lemma}[section]
\newtheorem{proposition}[lemma]{Proposition}

\theoremstyle{definition}

\newtheorem{definition}[lemma]{Definition}

\newcommand{\E}[1]{\left\langle {#1} \right\rangle}
\newcommand{\tr}{\operatorname{tr}}

\newcommand{\Aut}{\operatorname{Aut}}
\newcommand{\Ai}{\operatorname{Ai}}
\newcommand{\map}{\operatorname{map}}

\newcommand{\End}{\operatorname{End}}

\newcommand{\al}{\alpha}

\newcommand{\tsum}{{\textstyle\sum}}
\newcommand{\si}{{\boldsymbol\delta}}

\newcommand{\Mgn}{\overline{\mathcal{M}}_{g,n}}

\begin{document}

\title{Generating functions for intersection numbers
on moduli spaces of curves}
\author{Andrei Okounkov}
\date{June 2000}

\maketitle

\begin{abstract}
Using the connection between intersection theory on the  
Deligne-Mumford spaces $\Mgn$ and the edge scaling of
the GUE matrix model (see \cite{rp,OP}), we express the $n$-point functions
for the intersection numbers as $n$-dimensional error-function-type
integrals and also give a derivation of Witten's KdV equations
using the higher Fay identities of Adler, Shiota, and van Moerbeke. 
\end{abstract}

\section{Introduction}

\subsection{Overview}

\subsubsection{}

This paper is a continuation of \cite{rp} and is very 
closely related to \cite{OP}. It was observed
in \cite{rp} and conceptually explained in \cite{OP}
that the intersection theory on the moduli
space of curves is closely connected to the edge scaling of
the standard GUE matrix model (which can be called the 
edge-of-the-spectrum matrix model). Leaving the
detailed discussion of this phenomenon to \cite{OP}, we obtain here
some formal consequences of this connection, most importantly, an 
error-function-type integral
formula for certain generation functions for the intersection numbers
known as the $n$-point functions. We also show how to derive, using results
of Adler, Shiota, and van Moerbeke \cite{ASV}, the KdV
equations from this matrix model. 

\subsubsection{} 

Let $\Mgn$ be the stable compactification of the
moduli space of genus $g$ curves with $n$ marked points.
Let $\psi_i$ be the first Chern class of the line
bundle whose fiber over each pointed stable
curve is the cotangent line at the $i$th point.
We use the standard notation 
\begin{equation}\label{Et}
\E{\tau_{d_1} \dots \tau_{d_n}} = \int_{\Mgn} \psi_1^{d_1} \cdots
\psi_n^{d_n} \,, \quad \sum {d_i}=3g-3+n  \,,
\end{equation}
for the intersection numbers of these $\psi$-classes.

The celebrated conjecture of Witten \cite{W} says that  for the following
generating function
$$
\bF(t_0,t_1,\dots)= \sum_{(k_0,k_1,\dots)} 
\E{\tau_{0}^{k_0} \tau_{1}^{k_1} \dots } \prod \frac{t_i^{k_i}}{k_i!}\,,
$$
the exponential $\exp(\bF)$ is a $\tau$-function 
for the KdV hierarchy in variables
$$
T_{2i+1}=\frac{t_i}{(2i+1)!!}\,.
$$
This conjecture was inspired by matrix models of $2$-dimensional quantum gravity,
see for example \cite{D,DGZ} for a survey. 

The proof of Witten's conjecture given by Kontsevich
\cite{K} uses another matrix model interpretation of the
generating function $\bF$.

\subsubsection{}

The purpose of this paper is to exploit yet another random matrix 
connection \cite{rp} to evaluate a different generating 
function for \eqref{Et}, called the $n$-point function, in a closed form.

\begin{definition} We call the following generating function
$$
\cF(\ix_1,\dots,\ix_n) = \sum_{g=0}^\infty  \cF_g(\ix_1,\dots,\ix_n)\,,
$$
where 
$$
\cF_g(\ix_1,\dots,\ix_n) = 
\sum_{\sum d_i = 3g-3+n} \E{\tau_{d_1} \dots \tau_{d_n}} \,  \prod \ix_i^{d_i}\,,
$$
the $n$-point function.
\end{definition}

The 2 and 3-point function were computed by Dijgraaf and Zagier,
respectively. Several
specialization of the $3$-point function can be found in the paper
\cite{CP}. Our main result, see Theorem \ref{t1} below, is that
the n-point function $\cF$ is a certain specific multivariate error function
(and hence, in particular, a multivariate hypergeometric function
in the sense of Gelfand, Kapranov, and Zelevinsky, see Section \ref{GKZ} ). 

\subsubsection{}
We recall, referring the reader to \cite{rp,OP} for details,
that the $n$-point functions turn out to be connected to the asymptotics
of the following averages over $N\times N$ Hermitian matrices
\begin{equation}\label{EN}
\E{\prod_{k=1}^n \tr H^{2[x_k N^{2/3}]}}_{N}\,, \quad N\to\infty\,,
\end{equation}
with respect to the standard Gaussian measure. Here  $[\,\cdot\,]$ stands
for the integer part. The asymptotics of \eqref{EN} is equivalent
to the knowledge of the distribution of eigenvalues of a random matrix
near the edge of the Wigner semicircle. This distribution is well known to be 
described by
the Airy ensemble, see for example \cite{TW}. We also note that 
this Airy behavior near the
edge of the spectrum is very universal for random matrix
ensembles, see for example \cite{Sosh}. 

\subsubsection{} 

By a standard application of the Wick formula, the asymptotics \eqref{EN}
is equivalent to the asymptotics in the following combinatorial enumeration
problem.  Consider a surface
$\Sigma$ and a map on it with $n$ cells, that is, a way to glue $\Sigma$
out of $n$ polygons by identifying the edges of the polygons in pairs. 
The asymptotics of the number of ways this can be
done, as $n$ is fixed and the perimeters of the polygons go to infinity at certain relative rates,
is encoded in the asymptotics of \eqref{EN}. 
\begin{figure}[!hbt]
\centering
\scalebox{.7}{\includegraphics{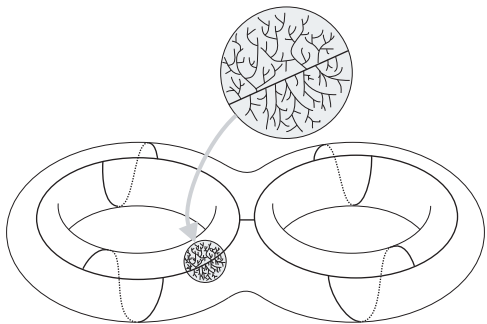}}
\caption{A typical map of large perimeter}
\label{f1}
\end{figure}
A typical map with 3 cells on a 
genus $g=2$ surface $\Sigma$ is shown in Figure \ref{f1}. Here one
portion of the map is magnified so that to show the dendriform pattern
that the map forms because the polygons are allowed to be glued to itself.
In fact, the overwhelming part of the perimeter is typically contained in 
these trees as the perimeter goes to infinity. The macroscopic data of
the map in Figure \ref{f1} is a trivalent graph embedded in $\Sigma$
which is the same data as in \cite{K}, see \cite{rp,OP} for details
of the connection of this approach to \cite{K}. 

The asymptotics in this map enumeration problem is described by
a certain function $\map_\Sigma(x_1,\dots,x_n)$, see Section 2.1.7 of
\cite{rp} for precise definitions and Section \ref{smap} below for
a formula for $\map_\Sigma$.
When the surface $\Sigma$ is connected of genus $g$, we write $\map_g$
instead of $\map_\Sigma$. Obviously, the asymptotics $\map_\Sigma$ is multiplicative
in connected components. 

There exist a number of other combinatorial asymptotics equivalent to
the asymptotics \eqref{EN}, such as the distribution of increasing
subsequences in a random permutation, see for example \cite{BDJ,O,BOO,J}.  

\subsubsection{}

Another product of our analysis is an alternative derivation of Witten's
KdV equations using the edge matrix model. 
It is based on the Fay identities methods developed
by Adler, Shiota, and van Moerbeke in \cite{ASV}.

\subsection{Formula for the $n$-point function}

\subsubsection{The function $\cE$}

The key ingredient in our formula for the $n$-point function 
will be the following function of $\ix_1,\dots,\ix_n$
\begin{multline}\label{forE}
\cE(\ix_1,\dots,\ix_n)=
\frac1{2^n\pi^{n/2}} 
\frac{\exp\left(\frac1{12} \sum \ix_i^3\right)
}{\prod \sqrt{\ix_i}} \times \\
\int_{s_i\ge 0} ds \,
\exp\left(-\sum_{i=1}^n \frac{(s_i-s_{i+1})^2}{4\ix_i} -\sum_{i=1}^n
\frac{s_i+s_{i+1}}2\,\ix_i 
\right) \,,
\end{multline}
where the integral is over $\R_{\ge 0}^n$ and $s_{n+1}=s_1$.

\begin{figure}[!hbt]
\centering
\scalebox{.7}{\includegraphics{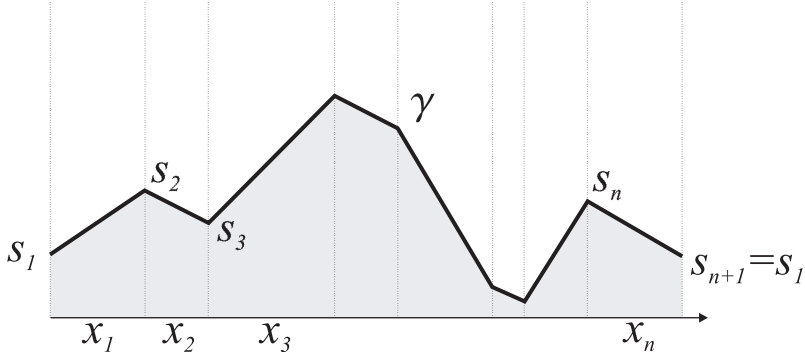}}
\caption{A piecewise linear function $\gamma$}
\label{f2}
\end{figure}

This integral admits a nice probabilistic interpretation, namely
$$
\cE(x) = \exp\left(\frac1{12} \sum x_i^3\right) \, \int_{\gamma\ge 0} 
e^{-\int \gamma} \, W(d\gamma)\,,
$$
where the integral is over the space of nonnegative
piecewise linear function of
the form shown in Figure \ref{f2}
subject to the bridge condition $s_{n+1}=s_{1}$.
The measure $W$ is the natural Gaussian measure: the increments
of $\gamma$ are independent normal variables with variance equal to 
twice the length of the interval.

\subsubsection{Main result}

The function $\cE$
is clearly invariant under a cyclic shift of the $\ix_i$'s.
We want to have something symmetric in the $\ix_i$'s, so
we set
$$
\cEt(\ix) =
\sum_{\sigma\in S(n)/(12\dots n)} \cE(\ix_{\sigma(1)},\dots,\ix_{\sigma(s)}) \,,
$$
where the summation is over coset representatives 
modulo the cyclic group generated by the permutation
$(12\dots n)$. 

Let $\Pi_n$ denote the set of all partitions $\al$ of the
set $\{1,\dots,n\}$ into disjoint union of subsets. For any
partition $\al\in\Pi_n$ with $\ell=\ell(\al)$ blocks, let $x_\al$ denote the  
vector of size $\ell$ formed by sums of $\ix_i$ over the blocks of $\al$.
For example, if
$$
\al=\{1,3\}\sqcup\{2\} \in \Pi_3
$$
then $\ell(\al)=2$ and $x_\al=(x_1+x_3,x_2)$.

By definition, we set 
$$
\cG(\ix_1,\dots,\ix_n) = \sum_{\al\in\Pi_n} (-1)^{\ell(\al)+1} \,  \cEt(\ix_\al)\,.
$$
For example, 
$$
\cG(\ix_1,\ix_2)= 
 \cE(\ix_1+\ix_2) -  \cE(\ix_1,\ix_2) \,,
$$ 
and, similarly, 
\begin{multline*}
\cG(\ix_1,\ix_2,\ix_3) =  \cE(\ix_1+\ix_2+\ix_3) - \cE(\ix_1+\ix_2, \ix_3) 
- \cE(\ix_1+\ix_3, \ix_2) \\
 - \cE(\ix_2+\ix_3, \ix_1)  +  \cE(\ix_1,\ix_2,\ix_3) +  \cE(\ix_1,\ix_3,\ix_2) \,.
\end{multline*}

In this notation, our result is
\begin{theorem}\label{t1}
\begin{equation}\label{main}
\cF(\ix_1,\dots,\ix_n) = \frac{(2\pi)^{n/2}}{\prod \ix_i^{1/2}} \, 
\cG\left(\frac{\ix}{2^{1/3}}\right) \,. 
\end{equation}
\end{theorem}

We remark that the right-hand side of \eqref{main} makes sense for
all $g\ge 0$ and $n\ge 1$, even for the pairs 
$$
(g,n)=(0,1),(0,2)
$$
for which the corresponding moduli space is problematic. The corresponding
terms, however, are not polynomial in $x$, see Section \ref{rmks}. 

\subsubsection{}\label{GKZ}

The integral \eqref{forE} fits inside a more general class of integrals
of the form
$$
I(Q)=\int_{\R_{\ge 0}^n} e^{Q(s)} ds\,, \quad \deg Q = 2 \,.
$$
Such a multivariate analog of the 
error function is, as a function of coefficients of the polynomial $Q$,
a multivariate hypergeometric function in the
sense of Gelfand, Kapranov, and Zelevinsky. In fact, this is true 
for any polynomial
$$
Q=\sum_{m=(m_1,\dots,m_n)} a_m  s^m \,,
$$
because the integral $I(Q)$ satisfies the obvious equations
$$
\left(\prod \frac{\partial}{\partial a_{\mu^{(i)}}} - 
\prod \frac{\partial}{\partial a_{\nu^{(i)}}}  \right)\, I(Q) = 0 \,,
\quad
\sum \mu^{(i)} = \sum \nu^{(i)}\,,
$$
and integration by parts gives
$$
I(Q)=-\int_{\R_{\ge 0}^n} e^{Q(s)} s_i \frac{\partial}{\partial s_i} Q \, ds\,,
$$
which is equivalent to the homogeneity equations 
$$
\left(\sum_{m} m_i \, a_m \frac{\partial}{\partial a_{m}} + 1\right)\, I(Q)= 0 \,,
$$
for $i=1,\dots,n$.

\section{Proof of the $n$-point function formula}

\subsection{}\label{smap}

We will assume that the reader is familiar  with Section 2 of \cite{rp}.
The main object we need from \cite{rp} is the function $\map_g(\ix)$ 
which describes the $N\to\infty$ asymptotics of the number of genus $g$ maps 
with $n$ cells with perimeters $\sim Nx_1,\dots , \sim N x_n$, respectively.

We have the following formula which is Theorem 3 in \cite{rp}
\begin{multline}\label{e1}
\int_{\R_{\ge 0}^n} e^{-(z,\ix)}\,\map_g(\ix)\,
\frac{d\ix}{\ix} =\\
\frac1{2^{|e(\Gamma)|/2-1}} \sum_{\Gamma\in \Gamma^{3}_{g,n}}
\frac{1}{|\Aut(\Gamma)|}  
\prod_{e\in e(\Gamma)} 
\frac1{\sqrt{z_{1,e}}+
\sqrt{z_{2,e}}}\,,
\end{multline}
where the summation is over all trivalent ribbon
graphs $\Gamma$ of genus $g$ with $n$ cells, the product is over
all edges $e$ of $\Gamma$, $z_{1,e}$ and $z_{2,e}$ are the $z_i$'s
corresponding to the two sides of edge $e$, and
$$
|e(\Gamma)| = 6g-6+3n
$$
is the number of edges of any graph $\Gamma\in \Gamma^{3}_{g,n}$. 

On the other hand, we have the following formula of Kontsevich
which is the unique boxed formula in \cite{K}
\begin{multline}\label{e2}
\sum_{\Gamma\in \Gamma^{3}_{g,n}}
\frac{1}{|\Aut(\Gamma)|}  
\prod_{e\in e(\Gamma)} 
\frac1{{z_{1,e}}+{z_{2,e}}} = \\
2^{-d-|e(\Gamma)|/3} 
\sum_{\sum d_i = 3g-3+n} \E{\tau_{d_1} \dots \tau_{d_n}} \, \prod_{i=1}^n
\frac{(2d_i)!}{d_i \, !} \, z_i^{-2d_i-1} \,.
\end{multline}
We refer the reader to \cite{OP} for a very detailed discussion and
proof of this formula. 

{F}rom \eqref{e1} and \eqref{e2}, using the formula 
$$
\int_0^\infty e^{-st}\, t^{m-\frac12} \, dt =
\Gamma\left(m+\frac12\right) \, s^{-m-\frac12} = \sqrt{\pi} \,
\frac{(2m)!}{2^{2m} \, m!} \,  s^{-m-\frac12} 
$$ 
and the fact that the function 
$\map_g$ is homogeneous of degree $3g-3+3n/2$,
we obtain the following:  

\begin{proposition}
\begin{equation}\label{Fmap}
\cF_g(\ix_1,\dots,\ix_n) = \frac{\pi^{n/2}}{2^g} \, 
\frac{\map_g(2\ix_1,\dots,2\ix_n)}{\prod \ix_i^{1/2}} \,. 
\end{equation}
\end{proposition}

\subsection{}

The function $\map_g$ can be obtained as the Laplace transform 
of a certain natural object. Introduce the Airy kernel, see
for example \cite{TW},
$$
\KK(x,y)= \frac{\Ai(x) \, \Ai'(y) - \Ai'(x) \Ai(y)}{x-y} \,,
$$ 
where $\Ai(x)$ is the classical Airy function. 

The kernel
$\KK$ differs by scaling of variables from the kernel $K$
used in \cite{rp}. It is natural to use $\KK$ when working with
$\map_g(2\ix)$ instead of $\map_g(\ix)$. We now introduce the
corresponding modifications of the functions from Section 2.1.8  of 
\cite{rp}. 

Let the function $\RR$ be the following Laplace transform
\begin{equation}\label{defR}
\RR(\xi_1,\dots,\xi_n)= \int_{\R^n} e^{(\xi,x)} \, \det \big[ \KK(x_i,x_j)\big)] \, dx\,.
\end{equation}
By definition, set
\begin{equation}\label{defH}
\HH(\ix_1,\dots,\ix_n) = \sum_{\al\in\Pi_n} \RR(\ix_\al)\,,
\end{equation}
where $\Pi_n$ denotes the set of all partitions of the
set $\{1,\dots,n\}$ into disjoint union of subsets, and 
 $\ix_\al$ is the  
vector formed by sums of $\ix_i$ over blocks of $\al$. 
For example,
\begin{multline*}
\HH(\ix_1,\ix_2,\ix_3) = \RR(\ix_1,\ix_2,\ix_3)+ \RR(\ix_1+\ix_2,\ix_3)+ \\
\RR(\ix_1+\ix_3,\ix_2) + \RR(\ix_2+\ix_3,\ix_1)+ \RR(\ix_1+\ix_2+\ix_3)\,.
\end{multline*}
Finally, set
\begin{equation}\label{defG}
\GG(\ix_1,\dots,\ix_n) = \sum_{S\subset\{1,\dots,n\}} \HH(\ix_i)_{i\in S} \,
\HH(\ix_i)_{i\notin S} \,,
\end{equation}
where the summation is over all subsets $S$ and $\HH(\ix_i)_{i\in S}$ denotes
the function $\HH$ in variables $\ix_i$, $i\in S$. For example
$$
\GG(\ix_1,\ix_2) = 2 \HH(\ix_1,\ix_2) + 2 \HH(\ix_1) \HH(\ix_2) \,.
$$

The function $\map_g$ has a natural extension to disconnected
surfaces $\Sigma$ by multiplicativity. 
We have the following formula, see the last formula in Section 2.1.8 of \cite{rp},
\begin{proposition} We have 
\begin{equation}\label{Su}
\GG(\ix_1,\dots,\ix_n) =  \sum_{\Sigma} 
\map_\Sigma(2\ix_1,\dots,2\ix_n) \,,
\end{equation}
where the summation is over all orientable surfaces $\Sigma$, 
including disconnected ones. 
\end{proposition}

For example, 
$$
\GG(\ix_1,\ix_2) = \sum_{g} \map_g(2\ix_1,2\ix_2) + 
\sum_{g_1,g_2} \map_{g_1} (2\ix_1)  \map_{g_2} (2\ix_2) \,.
$$

\subsection{}

\begin{definition} Let $Q(x_1,\dots,x_n)$, $n=1,2,\dots$,
be a sequence of function which are symmetric in their arguments.
We define 
\begin{align*} 
Q^\circ(x_1)&=Q(x_1)\,,\\
Q^\circ(x_1,x_2) &=Q(x_1,x_2) - Q(x_1) Q(x_2)\,, \\
Q^\circ(x_1,x_2,x_3) &= Q(x_1,x_2,x_3) -Q(x_1)Q(x_2,x_3) 
-Q(x_2)Q(x_1,x_3)  \\
&\qquad - Q(x_3)Q(x_1,x_2) + 2Q(x_1)Q(x_2)Q(x_3)\,,
\end{align*}
and so on, namely, terms of degree $d$ in $Q$ come with
the coefficient $(-1)^{d-1} (d-1)!$ which is the well known
M\"obius function for the partially ordered
set $\Pi_n$ of partitions of $\{1,\dots,n\}$. We call these 
functions $Q^\circ$ the \emph{connected} part of $Q$.  
\end{definition}

It is clear that the
connected part of the sum \eqref{Su} is the following 
$$
\Gc(\ix_1,\dots,\ix_n) = \sum_g \map_g(2\ix_1,\dots,2\ix_n) \,.
$$
It is clear from \eqref{Fmap} that we have
\begin{proposition}
$$
\cF(\ix_1,\dots,\ix_n) = \frac{2^{n/2-1}\, \pi^{n/2}}{\prod \ix_i^{1/2}} \, 
\Gc\left(\frac{\ix}{2^{1/3}}\right) \,. 
$$
\end{proposition}

Now it only remains to check that
\begin{equation}\label{cGGc}
\cG \overset{?}= \frac12 \, \Gc \,.
\end{equation}

\subsection{}

The sum $\GG$ can be interpreted in a diagrammatic way as follows. 
We take our variables $\{\ix_i\}$ and divide them first in two 
large groups
in all possible ways as in \eqref{defG}. Next we divide each group
into smaller  subgroups, which are the blocks of the partition $\al$ 
in \eqref{defH}. These blocks are then  grouped together
into the cycles of permutations which appear in  the determinant
in \eqref{defR}. 

The contribution of each such diagram to $\GG$
is a product over the connected pieces of the diagram. Therefore,
in the function $\Gc$ every disconnected diagram will cancel out
and $\Gc$ will be the sum over connected diagrams only. 

The corresponding summands are the following. First, in \eqref{defG}
every nontrivial subset $S$ leads to a disconnected diagram and
so the connected part of $\GG$ is just twice the connected
part of $\HH$. After that, in \eqref{defH}
we can take any partition $\al$, but then in \eqref{defR} we need
our permutation to have only one cycle. 

Denote the contribution of  one long cycle to \eqref{defR} by
\begin{equation}\label{defE}
\EE(\ix_1,\dots,\ix_n)=
 \int_{\R^n} e^{(\ix,z)} \, \prod_{i=1}^{n} \KK(z_i,z_{i+1}) \, dz \,,
\end{equation}
with the understanding that $x_{n+1}=x_1$. Note the invariance
of the  function
\eqref{defE} under the cyclic permutation  of variables.

The connected part of
$\RR$ is therefore
$$
\Rc(\ix) = (-1)^{n+1}
\sum_{\sigma\in S(n)/(12\dots n)} \EE(\ix_{\sigma(1)},\dots,\ix_{\sigma(s)}) \,,
$$
where the summation is over representatives of the right cosets of the
symmetric group modulo the cyclic group generated by the permutation
$(12\dots n)$. With this notation, we have the following 
\begin{proposition}\label{GcR}
$$
\Gc(\ix_1,\dots,\ix_n) = 2 \sum_{\al\in\Pi_n} \Rc(\ix_\al) \,.
$$
\end{proposition}

For example, 
$$
\Gc(\ix_1,\ix_2)= 2 \Rc(\ix_1,\ix_2) + 2 \Rc(\ix_1+\ix_2) =
2 \EE(\ix_1+\ix_2) - 2 \EE(\ix_1,\ix_2) \,.
$$ 
Similarly
\begin{multline*}
\Gc(\ix_1,\ix_2,\ix_3) = 2 \EE(\ix_1+\ix_2+\ix_3) -2 \EE(\ix_1+\ix_2, \ix_3) 
-2 \EE(\ix_1+\ix_3, \ix_2) \\
 -2 \EE(\ix_2+\ix_3, \ix_1)  + 2 \EE(\ix_1,\ix_2,\ix_3) + 2 \EE(\ix_1,\ix_3,\ix_2) \,.
\end{multline*}

\subsection{}

Now it only remains to prove that
$$
\cE \overset{?}= \EE \,.
$$

Consider the function \eqref{defE} more closely. We have the following
formula, see formula (4.5) in \cite{TW}, 
\begin{equation}\label{AA}
\KK(z,w)=\int_0^\infty \Ai(z+a) \, \Ai(w+a) \, da \,.
\end{equation}
Substituting this into \eqref{defE} and interchanging the order of
integration we obtain the integrals of the form considered in the 
following 

\begin{lemma} 
$$
\int_{-\infty}^\infty e^{\ix z} \Ai(z+a) \, \Ai(z+b) \, dz =
\frac{1}{2\sqrt{\pi \ix}} \, 
\exp\left(\frac{\ix^3}{12}-\frac{a+b}2\,\ix-\frac{(a-b)^2}{4\ix}\right)
$$
\end{lemma}

\begin{proof} Denote by $f$ the
function $f= \Ai(z+a) \, \Ai(z+b)$. The Airy function is
a solution of the Airy differential equation $\Ai''(z) = z \Ai(z)$. 
It follows that
$$
f_{zzzz}-(4z+2a+2b) f_{zz}-6 f_{z}+(a-b)^2 f =0 \,.
$$
This translates into a first order ODE for $g=\int e^{\ix z} f \, dz$. 
$$
g_\ix=\left(\frac{\ix^2}4-\frac{a+b}2-\frac1{2\ix}+\frac{(a-b)^2}{4\ix^2}\right) g \,, 
$$
from which it follows that
$$
g=\frac{h(a,b)}{\sqrt{\ix}} \, 
\exp\left(\frac{\ix^3}{12}-\frac{a+b}2\,\ix-\frac{(a-b)^2}{4\ix}\right)\,,
$$
where $h(a,b)$ is some function of $a$ and $b$. The obvious equations
$$
f_{aa} - (z+a) f =0\,, \quad f_{bb} - (z+b) f =0 
$$
translate into the equations
$$
g_{aa} - a g = g_{bb} - b g = g_\ix \,, 
$$ 
which imply that $h$ is constant. The fact that $h=\dfrac{1}{2\sqrt{\pi}}$ 
follows by taking the Laplace transform of \eqref{AA} with $z=w$ 
and comparing it with
the known formula for it, see Section 2.6 of \cite{rp}.  
\end{proof}

This lemma concludes the proof of the theorem. 

\subsection{Examples and remarks}

\subsubsection{}

In the simplest example $n=1$, there is no quadratic term in
the exponential and we obtain the formula
$$
\cE(\ix)= \frac1{2\sqrt{\pi}} 
\frac{\exp\left(\frac1{12} \ix^3\right)}{\sqrt\ix} \, \int_0^\infty e^{xy} \, dy 
=\frac1{2\sqrt{\pi}} 
\frac{\exp\left(\frac1{12} \ix^3\right)}{\ix^{3/2}}\,,
$$
which corresponds to the well-known result
$$
\E{\tau_{3g-2}} = \frac1{24^g\, g!}\,, \quad \cF(x)=\frac{\exp(x^3/24)}{x^2}\,. 
$$

\subsubsection{}
The case $n=2$ reduces to the ordinary error function
$$
\erf(x)=\frac2{\sqrt\pi}\int_0^x e^{-u^2} \, du \,,
$$
as follows. 
The following integral 
\begin{multline}\label{I2}
\iint_0^\infty \exp(-P(a+b)-Q(a-b)^2)\, da\, db = \\
\frac{\sqrt{\pi} }{2 P \sqrt{Q}} \, \exp\left(\dfrac{P^2}{4Q}\right) \, 
\left(1-\erf\left(\frac{P}{2\sqrt{Q}}\right)\right) 
\end{multline}
can be computed by introducing the variables $u=a+b$ and $v=a-b$ and
integrating first in $u$ and then in $v$. 

The integral in \eqref{forE} for $n=2$ becomes the integral \eqref{I2}
with the following choice of parameters
$$
P=\frac{\ix_1+\ix_2}{2} \,, \quad Q= \frac{\ix_1+\ix_2}{4\ix_1\ix_2} \,.
$$  
Therefore
$$
\cE(\ix_1,\ix_2)= \frac{1}{2\sqrt\pi} \frac{\exp\left((\ix_1+\ix_2)^3/12\right)}
{(\ix_1+\ix_2)^{3/2}} \, \left(1-\erf\left(\frac12\sqrt{\ix_1\ix_2(\ix_1+\ix_2)}\right)\right) 
$$
It follows that 
\begin{multline}\label{Gc2}
\cG(\ix_1,\ix_2) =  \cE(\ix_1+\ix_2) - \cE(\ix_1,\ix_2)  \\ 
= \frac{1}{2\sqrt\pi} \frac{\exp\left((\ix_1+\ix_2)^3/12\right)}
{(\ix_1+\ix_2)^{3/2}} \, \erf\left(\frac12\sqrt{\ix_1\ix_2(\ix_1+\ix_2)}\right) \,.
\end{multline}
Now we use the expansion
$$
e^{x^2/4} \, \erf(x/2) = \frac{1}{\sqrt\pi}
\sum_{k=0}^\infty \frac{k!}{(2k+1)!} \, x^{2k+1} \,,
$$
which can be proved, for example, by noting that the function $\erf(x/2)$
satisfies the equation $u_{xx} + \frac{x}2 \,  u_x = 0$ and hence the
function $e^{x^2/4} \, \erf(x/2)$ satisfies the equation 
$u_{xx} - \frac{x}2 \,  u_x - \frac12\, u = 0$. It follows that
\begin{equation}\label{F2}
\cF(\ix_1,\ix_2) = \frac1{\ix_1+\ix_2} 
\, \exp\left(\frac{\ix_1^3}{24} + \frac{\ix_1^3}{24}\right) 
\sum_{k=0}^\infty  \frac{k!}{(2k+1)!} \, \left(\frac12\,
 \ix_1\ix_2(\ix_1+\ix_2)
\right)^k \,.
\end{equation}
By the string equation
$$
\cF(\ix_1,\ix_2) = \frac{\cF(\ix_1,\ix_2,0)}{\ix_1+\ix_2}\,,
$$
and so again we find agreement with the known formula for $\cF(\ix_1,\ix_2,0)$,
see \cite{CP}. 

\subsubsection{}\label{rmks}

Observe that the only non-polynomial term in \eqref{F2}, which is
$\dfrac1{\ix_1+\ix_2}$  comes from
the exceptional case $(g,n)=(0,2)$. 

Also observe that the term $\cE(\ix_1+\ix_2)$ cancels out completely
in \eqref{Gc2} with a part of $\cE(\ix_1,\ix_2)$. This cancellation can be seen
a-priory as follows: all coefficient of $\cF$ are rational numbers
and the term $\cE(\ix_1+\ix_2)$ has the wrong power of $\pi$. This
argument works in general to identify terms making no
contribution to $\cG$. 

In fact, the formula \eqref{main} is full of cancellations of this
type. Several criteria for identifying irrelevant terms will be 
developed in the next section where we will be dealing with 
the KdV equations.

\section{KdV hierarchy}

\subsection{Strategy}\label{select}

\subsubsection{}

In this section we will give a derivation of the Witten's KdV equations 
which uses our current machinery and the Fay identities techniques
developed by Adler, Shiota, and van Moerbeke in \cite{ASV}. Namely,
we will prove that  the series expansion
of $\cF(x)$ about $x=0$ is  given
by coefficients of one specific KdV $\tau$-function. See \cite{K},
and also for example \cite{D},
for the exposition of how this was established originally
using Kontsevich's matrix model. 

\subsubsection{}

We will need some qualitative facts about the series expansion
of $\cG(x)$ about $x=0$, the first one being that 
all exponents in that expansion
are half-integers. Consequently, after sufficiently
many differentiations, every terms in this expansion
blows up as $x_i\to +0$ for any $i=1,\dots,n$. 
In other words, the coefficients of $\cG(x)$ can be 
determined from singularities of $\cG$ and its derivatives
as $x_i \to +0$.  

\subsubsection{}

Proposition \ref{GcR} and  \eqref{cGGc} express
the function $\cG$ in  in terms of the function $\RR$
which, by definition, is the following Laplace transform 
\begin{equation}\label{Lapl}
\RR(\xi_1,\dots,\xi_n)= \int_{\R^n} e^{(\xi,x)} \, \det \big[ \KK(x_i,x_j)\big]  \, dz\,.
\end{equation}
Hence the $x\to 0$ singularities of $\cG$ and its derivatives
are determined by the similar singularities of $\RR$. Clearly, 
the singularities of $\RR$ and its derivatives as $\xi_i \to +0$,
are determined by the $x\to -\infty$ asymptotics of the integrand in 
\eqref{Lapl}.  

The $x\to -\infty$ asymptotics of the integrand in \eqref{Lapl} can be
computed from the classical asymptotics of the Airy function, however
it is rather complicated and difficult to use in \eqref{Lapl} directly.
Fortunately, there are several qualitative criteria, a sort of 
selection rules, which help
 identify many terms in the asymptotics of $\det \big[ \KK(x_i,x_j)\big]$
as irrelevant.

\subsubsection{}\label{diagonal}

Observe that because all exponents in the expansion of $\cG(x)$ are 
half-integers, no information is lost by computing 
$\left(\frac{\partial}{\partial x_i} - \frac{\partial}{\partial x_j}\right)
\cG(x)$ instead of $\cG$. Indeed, the operator 
$$
\tfrac{\partial}{\partial x_i} - \tfrac{\partial}{\partial x_j} \in
\End_\C \left(\sqrt{x_1\cdots x_n}\, \C[x_1^{\pm1},\dots,x_n^{\pm1}]\right)
$$
has no kernel. On the other side of the Laplace transform \eqref{Lapl}, 
this means that we are free to multiply the integrand by any factor
of the form $(x_i-x_j)$.  

This is convenient for the following reason. 
The expression of $\cG$ in terms of $\RR$ involves terms which are
Laplace transforms of distributions supported on diagonals like
the $\cE(x_1+x_2)$ term in \eqref{Gc2}. Also, the asymptotics of
the integrand in \eqref{Lapl} contains denominators of the form
$\frac{1}{x_i-x_j}$
which require special handling resulting in 
certain contributions to the asymptotics from the diagonals.
Since these contributions disappear after multiplying by appropriate
factors of the form $(x_i-x_j)$, they always balance out in the final
answer like they did in \eqref{Gc2}. 

In other words, this
argument shows that if a term in the asymptotics of
$\det \big[ \KK(x_i,x_j)\big]$ becomes negligible after multiplication
by $(x_i-x_j)$, then it can be discarded.

We will call such negligible terms the \emph{diagonal} terms. For example,
Theorem \ref{t1} can be restated as saying that
\begin{equation*}
\cF(\ix_1,\dots,\ix_n) = \frac{(-1)^{n+1}(2\pi)^{n/2}}{\prod \ix_i^{1/2}} \, 
\cEt\left(\frac{\ix}{2^{1/3}}\right) + \, \textup{diagonal terms} \,. 
\end{equation*}
Another qualitative
way to eliminate the diagonal contributions is to observe that they
come with the wrong power of $\pi$, see Section \ref{rmks}. 

\subsubsection{}\label{osc}

Another category of negligible terms are the \emph{oscillating} terms in
the asymptotics of $\det \big[ \KK(x_i,x_j)\big]$. The $x\to-\infty$
asymptotics of the Airy function $\Ai(x)$ is oscillating and most
of the terms in the asymptotics of the integrand in \eqref{Lapl}
will have factors like $\exp\left(\frac{4}{3} i x_k^{3/2}\right)$ for some
$k=1,\dots,n$. We now observe that such a term will never blow
up as $\xi_k\to +0$, even if we multiply it by an 
arbitrary large power of  $x_k$. Indeed, we can assume that we already
got rid of all denominators as 
 explained in \ref{diagonal} above, and then it is enough to show that
$$
\int_c^\infty e^{-\xi x + i x^{3/2}} x^a \, dx = O(1)\,, \quad \xi\to + 0 \,,
$$
where $c>0$ and $a$ is arbitrary. After a change of variables, it
becomes equivalent to
$$
\int_c^\infty e^{-\xi x^{2/3}}  e^{i x} x^a \, dx = O(1)\,, \quad \xi\to + 0 \,,
$$
for some other constants $a$ and $c$. Now, if we integrate by parts 
integrating $ e^{i x}$ and differentiating the rest, we can decrease $a$
so that to make the integral with $\xi=0$ absolutely converging, thus proving
our assertion. 

\subsubsection{}\label{posdeg}
Finally, we are interested in terms in $\cG(x)$ of
 positive degree in all $x_i$.
Recall that all terms of $\cG(x)$ have positive degree in 
all $x_i$'s with the exception of the $g=0$ terms for $n=1,2$.

\subsection{Asymptotics of the Airy kernel}

\subsubsection{}

Introduce the following functions
$$
a(x)=\sum_{k=0}^\infty a_k x^{-3k} \,, \quad \ad(x)=\sum_{k=0}^\infty \ad_k x^{-3k+1}\,,
$$
where $a_0=\ad_0=1$ and 
$$
a_k = (-i)^k \, \frac{(6k-1)!!}{72^k (2k)!} \,, \quad \ad_k= \frac{1+6k}{1-6k}\, a_k \,.
$$ 
Also set, by definition,
$$
A(x)=\exp\left(i\tfrac23x^3\right) a(x)\,, \quad \Ad(x)=\exp\left(i\tfrac23x^3\right) \ad(x) \,.
$$ 
and
$$
\af(x,y)=\frac{a(x)\, \ad(y)- a(y)\,  \ad(x)}{x-y} \,, \quad
\A(x,y)=\frac{A(x)\, \Ad(y)- A(y)\,  \Ad(x)}{x-y} \,, 
$$ 
We have (see e.g.\ \cite{Hand}) the following $x\to+\infty$ asymptotics 
\begin{alignat*}{2}
\Ai(-x)&\sim \frac1{2\sqrt{\pi} i x^{1/4}} &&\Big[e^{\pi i/4} A\left(x^{1/2}\right)
- e^{-\pi i/4} A\left(-x^{1/2}\right) \Big] \,, \\
\Ai'(-x)&\sim - \frac1{2\sqrt{\pi} x^{1/4}} &&\Big[e^{\pi i/4} \Ad\left(x^{1/2}\right)
- e^{-\pi i/4} \Ad\left(-x^{1/2}\right) \Big] \,. 
\end{alignat*}
This asymptotics remains valid for complex $x$ such that $|\arg x|<\frac23\pi$. 

\subsubsection{}
It follows that away from the diagonal $x=y$ we have
$$
\KK(-x,-y)\sim\frac{1}{4\pi i \, x^{1/4} y^{1/4}} 
\sum_{\e,\e'=\pm 1}
i^{(\e+\e')/2} \,
\frac{\A\left(\e \sqrt x,\e' \sqrt y\right)}{\e \sqrt x+\e'\sqrt y}\,,
$$
where the sum is over $4$ possible combinations of the signs $\e$ and $\e'$. 
Because the function $\KK(-x,-y)$ is analytic at $x=y$ and
its asymptotics remains valid in the complex domain,
this asymptotics can be extended, using contour integrals,
to the diagonal $x=y$.

By multinearity of the determinant, we have
\begin{multline}\label{asK}
\det \big[\KK(-x_i,-y_j)\big]_{1\le i,j\le n} \sim  \\
\frac{1}{(4\pi i)^n  \prod_i x_i^{1/4} y_i^{1/4}}
\sum_{\e,\e'\in \{\pm 1\}^n}  
\det \left[\frac{\A\left(\e_i \sqrt{x_i},\, \e'_j  \sqrt{y_j}\right)}
{\e_i \sqrt{x_i}+ \e'_j  \sqrt{y_j}}
\right]\,,
\end{multline}
where the summation is over $4^n$ choices of signs vectors $\e$ and $\e'$. 
Observe that by definition of $\A$ we have 
$$
\det \left[\frac{\A\left(x_i, y_j \right)}{x_i + y_j}
\right] = 
\exp\left(\tfrac23 i\tsum (x_i^3 + y_i^{3}) \right)
\det \left[\frac{\af\left(x_i, y_j \right)}{x_i + y_j}
\right]\,.
$$
This last determinant can be simplified using the higher Fay identities \cite{ASV}
which will be discussed in the following section.

\subsection{Fay identities and $\tau$-functions}

\subsubsection{} 

In this section we collect, for the reader's convenience, some
background material about $\tau$-functions and Fay identities \cite{ASV}.  
We will slightly deviate from the 
notational conventions of the book \cite{Ka} by V.~Kac.
For our purposes, it will be convenient
to use the notations of the Appendix to \cite{O}. 

\subsubsection{} 

Given a matrix $M\in GL(\infty)$, we will denote by $\E{M}$ the
charge $0$ vacuum matrix element
$$
\E{M}=(M v_\emptyset,v_\emptyset)   
$$ 
for the action of $GL(\infty)$ in the infinite wedge space. For
any $M$, the function
$$
\tau_M(t) = \E{\Gamma_+(t) \, M}\,,
$$ 
of the variables $t=(t_1,t_2,\dots)$, is a $\tau$-function for the
KP hierarchy. In our case, we will additionally assume that $M$
is upper triangular, which implies that $M^* v_\emptyset=v_\emptyset$
and hence 
$$
\tau_M(0)=1 \,.
$$ 

\subsubsection{} 

Consider the following matrix element
$$
\Psi=\E{\prod_{i=1}^n \psi(z_i) \, \psi^*(w_i) \, M}\,,
$$
where $\psi(z)$ and $\psi^*(w)$ are the standard generating functions 
for the fermionic operators. On the one hand we have (see Corollary
14.10 in \cite{Ka} reproduced in (A.14) in \cite{O})
$$
 \psi(z) \, \psi^*(w) = \frac{\sqrt{zw}}{z-w} \, \Gamma_-(\{z\}-\{w\}) \, 
\Gamma_+(\{w^{-1}\}-\{z^{-1}\})\,,
$$
where
$$
\{z\}=\left(z,\frac{z^2}{2},\frac{z^3}{3},\dots\right) \,.
$$ 
{}From this and the commutation rule
$$
\Gamma_+(t)\, \Gamma_-(s) = e^{\sum k t_k s_k} \,  \Gamma_-(s)\, \Gamma_+(t)
$$
it follows that
\begin{equation}\label{fay1}
\Psi=\Pi_0 \,  \frac{\Delta(z) \, \Delta(-w)}{\prod (z_i-w_j)} \,\, 
\tau\left(\tsum \{w_i^{-1}\}-\tsum \{z_i^{-1}\}
\right) \,,
\end{equation}
where $\Delta$ denotes the Vandermonde determinant and $\Pi_0=\prod_i \sqrt{z_i w_i}$.

\subsubsection{} 

On the other hand, consider the operators $\psi_M(z)=M^{-1} \psi(z) M$,
and similarly $\psi^*_M(w)$. Using the commutation rules and
$M^* v_\emptyset=v_\emptyset$, we compute
\begin{multline}\label{fay2}
\Psi=\E{\psi_M(z_1)\cdots\psi_M(z_n) \, \psi^*_M(w_n) \cdots \psi^*_M(w_1)}=\\
\det\big[\E{\psi_M(z_i) \, \psi^*_M(w_j)}\big] = 
\Pi_0 \, \det\left[ 
\frac{\tau\left(\{w_i^{-1}\}-\{z_j^{-1}\}\right)}{z_i-w_j}\right] \,.
\end{multline}
Here the second equality is based on Wick's theorem, or equivalently,
follows from the observation that, acting on the vacuum, the operators
$\psi^*_M(w_j)$ remove vectors $\underline{k}$,  where $k\in\{-\frac12,-\frac32,\dots\}$,
and then the operators $\psi_M(z_i)$ have to put all these removed
vectors back,
in all possible orders. 

Combining \eqref{fay1} with \eqref{fay2} and interchanging the roles
of $z$ and $w$, we obtain the identity
\begin{multline}\label{fay3}
\det\left[ 
\frac{\tau\left(\{z_i^{-1}\}-\{w_j^{-1}\}\right)}{z_i-w_j}\right] = \\ 
\frac{\Delta(z) \Delta(-w)}{\prod (z_i-w_j)} \,\, 
\tau\left(\tsum \{z_i^{-1}\}-\tsum \{w_i^{-1}\}
\right)\,,
\end{multline}
which is the form of the higher Fay identities \cite{ASV} that we will
need here. 

\subsubsection{} 

{}From now on, we will be interested in one particular $\tau$-function, see
e.g.\ \cite{KS}, corresponding to the matrix $M_A$ which acts as follows
$$
M_A \, \underline{k} = \begin{cases}
\sum a_i \, \underline{k+3i}\,, & \textup{$k+1/2$ is even}\,,\\
\sum \ad_i \, \underline{k+3i}\,, & \textup{$k+1/2$ is odd}\,.
\end{cases}
$$
That is, the matrix $M$ looks like this: 
$$
M_A=
\begin{bmatrix}
\ddots&&& a_1 &&&  \ad_2\\
& a_0 &&& \ad_1 &&& a_2\\
&&\ad_0&&& a_1 \\
&&& a_0 &&& \ad_1 \\ 
&&&& \ad_0 &&& a_1\\
&&&&& a_0 \\
&&&&&& \ad_0 \\
&&&&&&& \ddots
\end{bmatrix}
$$
where the empty spaces represent zeros. This matrix commutes
with the bosonic operators $\alpha_k$ for $k$ even and hence the
$\tau$-function
$$
\tau_A = \tau_{M_A}
$$
does not depend on $t_k$ with $k$ even, meaning that it is a 
$\tau$-function for the KdV hierarchy. It also follows that
\begin{equation}\label{fay4}
\tau_A(t-\{w\})=\tau_A(t+\{-w\}) 
\end{equation}
for any $w$. 

\subsubsection{} 

{}It follows from definitions that
$$
\tau_A(\{x^{-1}\}+\{y^{-1}\})=\sum_{i,j=0}^{\infty}
a_i\, \ad_j \, \frac{x^{-3i-1} y^{-3j} - y^{-3i-1} x^{-3j}}{x^{-1}-y^{-1}} =
- \af(x,y) \,,
$$ 
which by \eqref{fay4} and the Fay identity \eqref{fay3} implies that
\begin{equation}\label{fay5}
\det \left[\frac{\af\left(x_i, y_j \right)}{x_i + y_j}
\right] = (-1)^n \frac{\Delta(x)\, \Delta(y)}{\prod (x_i+y_j)} \, \, 
\tau_A\left(\tsum \{x_i^{-1}\}+\tsum \{y_i^{-1}\}
\right)\,.
\end{equation}

\subsection{Asymptotics and KdV equations}

\subsubsection{}

Combining the formulas \eqref{asK} and \eqref{fay5}, we obtain 
the following asymptotics
\begin{multline}\label{asmptK}
\det \left[\KK\left(-x_i^2,-y_j^2\right)\right]_{1\le i,j\le n} \sim  
\\
\frac1{(-4\pi i)^n  \prod_i \sqrt{x_i y_i}} 
\sum_{\e,\e'\in \{\pm 1\}^n}  e^{\frac23 i\sum (\e_i x_i^{3} + \e'_i y_i^{3})} \times \\
\frac{\Delta(\e x)\, \Delta(\e' y)}
{\prod \left(\e_i x_i+\e'_j y_j\right)} 
\,
\tau_A\left(\tsum \left\{\frac{\e_i}{x_i}\right\}+
\tsum \left\{\frac{\e'_i}{y_i}\right\}
\right)\,, 
\end{multline}
where $\e x$ denotes the term-wise product of two vectors.

Since the left-hand side in \eqref{asmptK} is analytic and \eqref{asmptK}
is valid in the complex domain, it can be extended to the asymptotics 
 of the integrand in \eqref{Lapl} 
by letting $y_i \to x_i$ in \eqref{asmptK}. 

The full form of the $y_i \to x_i$ limit of \eqref{asmptK} is rather
messy, but fortunately, we can concentrate only on a small fraction
of the terms that survive all selection rules discussed in
Section \ref{select}. 

\subsubsection{}

Most importantly, we get nonoscillating terms only if $\e' = - \e$
which cuts the number of summands from $4^n$ down to
$2^n$. These $2^n$ terms can be conveniently organized using the 
following operators $\nabla_i$. 

Let $f(x_i,y_i,\dots)$ be  a function
of $x_i$, $y_i$ and some other variables which is 
supersymmetric in $x_i$ and $y_i$ in the sense that
its restriction to the diagonal 
$$
 \si_i f = f\big|_{y_i=x_i}
$$
does not depend on $x_i$, that is,
\begin{equation}\label{ss}
\frac{\partial}{\partial x_i} \si_i  f = 0 \,.
\end{equation}
Introduce the following operator
$$
\big[\nabla_1 f\big](x_1,\dots) = 
\lim_{y_1\to x_1}
 \frac12 \frac{f(x_1,y_1,\dots)-f(-x_1,-y_1,\dots)}{x_1-y_1} \,.
$$ 
Note that the numerator here vanishes on the diagonal $y_1=x_1$,
so this limit is well defined. 

With this notation, we see that
\begin{equation}\label{asmptK2}
\det \left[\KK\left(-x_i^2,-x_j^2\right)\right]_{1\le i,j\le n} \sim 
\frac1{(-2\pi i)^n  \prod_i x_i} \, \nabla_1 \nabla_2  \cdots \nabla_n \, \Phi + \dots
\end{equation}
where 
\begin{align}\label{ph}
\Phi&= \prod_i \phi_i \prod_{i\ne j} \phi_{ij} 
\,
\tau_A\left(\tsum \left\{\frac1{x_i}\right\}-
\tsum \left\{\frac{1}{y_i}\right\}
\right)\,, \\
\phi_i &= \exp\left(\tfrac23 i (x_i^{3}-y_i^{3})\right)\,, \quad
\phi_{ij} = \frac{(x_i-x_j)(y_j-y_i)}{(x_i-y_j)(x_j-y_i)}\,, \label{ph2}
\end{align}
and dots in \eqref{asmptK2} stand for oscillating terms. Observe that 
all factors in \eqref{ph}, \eqref{ph2} have the 
supersymmetry \eqref{ss}. 

\subsubsection{} The operators $\nabla_i$ commute and satisfy the
following Leibnitz-like rule
\begin{equation}\label{Leib}
\nabla_i (f\cdot g) = \nabla_i (f) \cdot \si_i(g) +
\si_i(f) \cdot   \nabla_i (g) \,.
\end{equation}
Let us examine the effect of applying $\nabla_i$ to \eqref{ph}. 

Denote by $\tau_\mu$ the coefficients in the expansion of the
$\tau_A$
\begin{equation}\label{taucoeff}
\tau_A = \sum_\mu \frac{\tau_\mu \, t_\mu}{|\Aut \mu|}
\end{equation}
where the summation is over all partitions $\mu$, $|\Aut \mu|$ is 
the product of factorials of multiplicities of parts in $\mu$, 
$t_\mu=\prod t_{\mu_i}$. The variables $t_k$ are 
specialized in \eqref{ph} in the following way
\begin{equation}\label{txy}
t_k = \frac1k\sum x_i^{-k}  -  \frac1k\sum y_i^{-k} \,.
\end{equation}
Recall that $\tau_A$ depends only on $t_k$ with $k$ odd. 
It is clear that
\begin{equation}\label{nabt}
\nabla_i \, t_k = - x_i^{-k-1}\,, \quad \textup{$k$ odd}\,,    
\end{equation}
and also that $\si_i$ simply removes $x_i$ and $y_i$ from the
sum \eqref{txy}. 

Similarly, it is clear that
\begin{equation}\label{nph}
\nabla_i \, \phi_i  = 2i x_i^2 \,,
\end{equation}
Finally, for the last type of factors in \eqref{ph} we have
$$
\si_i \, \phi_{ij}= \si_j \, \phi_{ij} = 1 
$$
which, in particular implies that
$$
\nabla_i\, \si_j  \, \phi_{ij} = \nabla_j \,\si_i  \, \phi_{ij} =0 \,.
$$
It follows that the $\nabla_i$'s have to be applied to $\phi_{ij}$'s in pairs
to get a nonzero result
\begin{equation}\label{nnph}
\nabla_i \nabla_j \, \phi_{ij} = \frac{x_i^2 + x_j^2}{(x_i^2-x_j^2)^2} \,.
\end{equation}

\subsubsection{}

The factor corresponding to \eqref{nnph} in the asymptotics of
$\det \big[\KK(-x_i,-x_j)\big]$ will be
\begin{equation}\label{nnph2}
\frac1{\sqrt{x_i x_j}}\frac{x_i + x_j}{(x_i-x_j)^2} \,, 
\end{equation}
where the first term comes from the prefactor in \eqref{asmptK2}.  
This has a second order pole on $x_i=x_j$, which in the 
full asymptotics of \eqref{asmptK2} cancels out with poles
of oscillating terms. This singularity of \eqref{nnph2} is immaterial 
and can be removed as explained in Section \ref{diagonal}. What
is important about \eqref{nnph2} is that is has degree $-2$ in
$x_i$ and $x_j$ and, hence, whatever contribution it makes to the
asymptotics of \eqref{Lapl}, it will be in degree $0$ in $\xi_i$ and $\xi_j$.
Since we are interested in terms of strictly positive degree in
all variables, see Section \ref{posdeg}, the terms containing \eqref{nnph2}
are negligible. By the same token, any terms containing \eqref{nph} can be
also discarded. 

It follows that relevant terms in the asymptotics \eqref{asmptK2} are
obtained by applying all operators $\nabla_i$ to the $\tau$-function.
{} From \eqref{nabt} and the Leibnitz rule \eqref{Leib} we obtain 
$$
\nabla_1 \cdots \nabla_n \, t_\mu =
\begin{cases}
(-1)^n |\Aut\mu| \, m_\mu\,, & \ell(\mu)=n \,, \\
0 & \textup{otherwise} \,,
\end{cases}
$$ 
where $m_\mu$ denotes the monomial symmetric function
$$
m_\mu = \frac1{|\Aut\mu|} \, \sum_{\sigma\in S(n)} \prod x_{\sigma(i)}^{\mu_i} \,.
$$ 
Therefore, modulo irrelevant terms, the asymptotics of the integrand in
\eqref{Lapl} is 
$$
\det \big[\KK(-x_i,-x_j)\big] \sim \frac{1}{(-2\pi i)^n \prod x_i} \sum_{\ell(\mu)=n}
\tau_\mu \, m_\mu(x^{-1/2}) + \dots \,. 
$$
Recall that only partitions $\mu$ with odd parts enter this sum.

\subsubsection{} 

Now it remains to use the formulas 
$$
\int_c^\infty e^{-\xi x} x^{a-1} \, dx \sim \frac{\Gamma(a)}{\xi^a} \,, \quad \xi\to+0\,,
$$
and
$$
\Gamma\left(-\frac{2k+1}{2}\right) = \frac{(-2)^{k+1} \sqrt{\pi}}{(2k+1)!!} 
$$
to obtain the expansion
\begin{align}\label{RRexp} 
\RR(\xi) &= \frac{\prod\xi^{1/2}}{i^n \,\pi^{n/2}} \sum_{m_1,\dots,m_n} 
\tau_{2m+1} \, \prod\frac{(-2\xi_i)^{m_i} }{(2m_i+1)!!} +\dots  
\end{align}
where dots stand for irrelevant terms, that is, for diagonal terms
 and terms of negative degree in the $\xi_i$'s, and
$\tau_{2m+1}$ is the coefficient in \eqref{taucoeff} corresponding to
 the nonincreasing rearrangement of
the numbers $2m_i+1$. 

By \eqref{cGGc} and Proposition \ref{GcR} we have
$$
\cG(x) = \sum_{\al\in\Pi_n} \Rc(\ix_\al)  = \Rc(x)+\dots\,,
$$
where $\Rc$ is the connected part of $\RR$ and dots stand
for  diagonal terms. Now \eqref{RRexp} and Theorem \ref{t1}
imply that
$$
\cF(x)=\left(\sum_{m_1,\dots,m_n} 
\wtt_{2m+1} \, \prod\frac{x_i^{m_i} }{(2m_i+1)!!} \right)^\circ + \dots\,,
$$
where $\wtt$ denotes the rescaled $\tau$-function
\begin{equation}\label{resc}
\wtt_\mu =  (-i)^{|\mu|} 2^{|\mu|/3} \, \tau_\mu \,,
\end{equation}
circle stands for the connected part, and dots stand for terms of negative
degree in $x$, that is, for the $g=0$ terms in the  case $n=1,2$.
This concludes the proof of the KdV equations.

Note that the rescaling \eqref{resc} is equivalent to the following
rescaling
$$
a_k \mapsto \frac{(6k-1)!!}{36^k (2k)!}\,, 
$$
keeping the relation $ \ad_k= \frac{1+6k}{1-6k}\, a_k$\,.

\end{document}